\documentstyle[amsfonts,amssymb,amsmath,12pt,latexsym]{article}

\newcommand{\sech}{\rm sech}
\newcommand{\dsp}{\displaystyle}
\newcommand{\pa}{\partial}
\newcommand{\de}{\delta}
\newcommand{\la}{\lambda}

\newcommand{\De}{\Delta}

\newcommand{\R}{\mathbb{R}}
\newcommand{\mb}{\mathbb{}}

\makeindex             

\newcommand{\beqst}{\begin{eqnarray*}}
\newcommand{\eeqst}{\end{eqnarray*}}
\newtheorem{theorem}{Theorem}
\newtheorem{remark}{Remark}

\newtheorem{lemma}{Lemma}

\newtheorem{definition}{Definition}
\date{}
\begin{document}
\title{Small Data    Wave Maps   in
Cyclic Spacetime}
\author{Karen Yagdjian\footnote{karen.yagdjian@utrgv.edu}, Anahit Galstian\footnote{anahit.galstyan@utrgv.edu} and Nathalie M. Luna-Rivera\footnote{nathalie.lunarivera01@utrgv.edu}\\
{}\\
 \small{School of Mathematical and Statistical Sciences,}\\
   \small{University of Texas RGV,   Edinburg, TX 78539, U.S.A.}
}

\maketitle

\begin{center}
Dedicated to Michael Reissig on his 60th birthday
\end{center}

\abstract{ 
    We study the initial value problem for the wave maps   defined on the cyclic spacetime with the target Riemannian manifold that is responsive (see definition of the self coherence structure) to the parametric resonance phenomena. In particular, for arbitrary  small  and smooth initial data we construct blowing up solutions of the wave map if the metric of the base manifold is periodic in time.}

\section{Introduction}

In this note we study   a wave map
\[
\phi :\, (L, g_{\mu \nu}) \longrightarrow (M, h_{ab})\,,
\]
where $L$ is an $n+1$-dimensional Lorentzian manifold
and the target $M$ is a $m$-dimensional Riemannian manifold. The map $\phi  $ is a {\it wave map} if it is a stationary point for the Lagrangian functional
\[
{\cal L}[\phi ]=\int_L\frac{1}{2} g^{\mu \nu }(x ) h_{ab} (\phi ) \nabla_\mu \phi ^a \nabla_\nu  \phi ^b  \, d\mu _g\,.
\]
The Lagrangian is written in local coordinates on the target, for which the notation $\phi^a =\phi^a (x^\mu)$  is used. We denote by $d\mu _g$ the measure with respect to the metric $g^{\mu \nu } $ on the spacetime.
Here the convention to  write $g^{\mu \nu }(x )=(g_{\mu \nu }(x ))^{-1} $ and $h^{ab} (\phi )=(h_{ab} (\phi ))^{-1} $ for the inverse of two metric tensors is used. These tensors are used  also in raising indexes. 
A stationary point for the Lagrangian functional  implies the following system of equations
\begin{eqnarray*}
\square u ^b - \Gamma ^b_{cd}(u)g^{\mu \nu }(x )   \nabla_\mu      u ^c \nabla_\nu       u ^d   =0\,,
\end{eqnarray*}
where $ \square$ is the d'Alembert (or wave) operator
\[
 \square := - \nabla _\mu \nabla ^\mu
\]
and $  \Gamma ^b_{cd}$ are the Christoffel symbols on the target manifold $(M,h)$ defined as 
\[
\Gamma_{j,k}^i (u) := \frac{1}{2}\sum_{l=1}^m h^{il} \left(\frac{\pa}{\pa u^j}h_{kl} + \frac{\pa}{\pa u^k}h_{jl} - \frac{\pa}{\pa u^l}h_{kj}\right)\,.
\] 
For the Minkowski spacetime ${\mathbb R}^{1+n}$ to a Riemanian manifold $M$ wave map satisfies the system of equations
\begin{equation}
\label{generalfns}
    \square u^i + \sum_{j,k=1}^m\Gamma_{j,k}^i (u) \left(\dot{u}^j\dot{u}^k- \nabla u^j \cdot \nabla u^k \right)=0\,,\quad i=1,\ldots m,
\end{equation}
where $\square=\partial^2/\partial t^2- \Delta $ and  $\Delta$ is the Laplacian in $L$. 
Here   $\dot{u}$ denotes the partial derivative with respect to time, and $\nabla$ denotes the gradient in $x$.

\bigskip

 For equation (\ref{generalfns})   consider the Cauchy problem with the initial conditions
 \begin{eqnarray}
\label{IC}
   &  &
   u^i(0,x) = u_0^i(x), \quad u_t^i(0,x) = u_1^i (x)\,, \quad i=1,\ldots m, \quad x \in {\mathbb R}^n\,.
\end{eqnarray}
It is known (see, e.g., Theorem 6.4.11 \cite{Hormander}) the following {\sl local existence result}:    if $\Gamma_{j,k}^i (u) $ are $C^\infty$ functions and
$u_0^i(x) \in H^{s+1} ({\mathbb R}^{n}) $ and $u_1^i(x) \in H^{s} ({\mathbb R}^{n}) $ for some integer $s> (n+2)/2$ then the  problem (\ref{generalfns})-(\ref{IC}) has for some $T>0$ a solution
$
u \in C^2 ([0,T]\times {\mathbb R}^n)$.
\smallskip

For the wave map from the  Minkowski spacetime ${\mathbb R}^{1+n}$, $n \geq 4$, to a Riemanian manifold $M$ the global in time existence of the small data solution can be derived  from Theorem 6.5.2~\cite{Hormander}. 
Klainerman and Machedon \cite{K-Machedon}  proved that the Cauchy problem for (\ref{generalfns}) is   locally  in time   well-posed in the Sobolev space $H^s({\mathbb R}^{1+n}) $ for any $s>n/2 $ if $\Gamma_{j,k}^i (u) $ are analytic and $n=3$.  Klainerman and Selberg~\cite{K-Selberg}  extended this result to $n \geq 2$.
\smallskip

Sideris~\cite{Sideris} considered wave maps (\ref{generalfns}) on the  Minkowski spacetime,  where $\Gamma_{j,k}^i (u) $ are smooth functions on ${\mathbb R}^m$ with the property  
\begin{equation}
\label{Sideris}
\Gamma_{j,k}^i (u^1,0,\cdots,0)=0 \quad \mbox{\rm for all}  \quad u^1 \in {\mathbb R},\,\, 1\leq i,j,k\leq m\,.
\end{equation}
Since the nonlinearities in (\ref{generalfns}) are cubic, small amplitude solutions are known to exist (see, e.g., \cite{Hormander}). In \cite{Sideris} the component $u^1$ need not to be small.
\smallskip

Georgiev and Schirmer in \cite{Georgiev}  generalized the spacetime estimates obtained by Klainerman and Machedon to wave equations on manifolds with nonconstant metric.  They applied these estimates to the question of global existence of low-regularity solution for small data of nonlinear wave equations on Minkowski space
${\mathbb R}^{1+3}$ satisfying the null condition. 
The null forms are expressions of the form
$g^{\mu \nu} \nabla_ \mu  u \nabla_ \nu  v$ or $\nabla_ \mu  u \nabla_ \nu  v-\nabla_ \nu  u \nabla_ \mu  v$, where $u,v$ are the functions on $ L$.    These estimates were  then applied on the Einstein cylinder (after  Penriose  compactification) to prove that if  $\left(u(0), u_{t}(0)\right) \in H^{2,1}({\mathbb R}^3) \times H^{1,2}({\mathbb R}^3)$  is sufficiently small, then a semilinear wave equations $\left(\partial_{t}^{2}-\Delta\right) u=F\left(u, \nabla u , u_{t}\right)$    with    $F $  satisfying  the null condition has a global  solution. 
\smallskip

In connection with low dimension $n$  we recall conjecture of Klainerman that states: {\it Let $({\mathbb H}^2,h)$ be the standard hyperbolic plane. Then classical wave maps originating on ${\mathbb R}^{2+1} $ exist for  arbitrary smooth initial data.}
\smallskip

The   answer to the Klainerman's conjecture  as well as the scattering result for the wave map  are given by Krieger and Schlag in \cite{Krieger,Krieger-Schlag}.
In particular, it is proved in \cite{Krieger-Schlag}
that if $M$ is a hyperbolic Riemann surface, and initial data $ (u(0),\partial_t u(0)) \,:\, S_0 \longrightarrow M\times TM$
are smooth and $u(0)=const$, $ \partial_t u(0)  =0$  outside of some compact set, then the wave map evolution $u $
of these data as a map ${\mathbb R}^{2+1} \longrightarrow M$ exists globally as a smooth function.
\smallskip

In \cite{Nishitani-Yagdjian} the stability of the last result under perturbation of the metric $g $ in $L$, that is,   in the perturbed Minkowski spacetime is investigated.
More exactly,    by Nishitani and Yagdjian \cite{Nishitani-Yagdjian}    considered the case of the  Riemannian manifold   $(M,h)$,  which  belongs to one-parameter family of manifolds containing
the Euclidean half-space and  the Poincar{\'e} upper half-plane model  $({\mathbb H}^2,h)$.
In fact, that family consists of the Riemannian manifolds, which are
 the half-plane $ \{(u^1,u^2) \in {\mathbb R}^2 \,|\, u^2>0 \}$ equipped with the metric   $\displaystyle h_{ij}du^i du^j = \frac{1}{(u^2)^l}\left( (du^1)^2+(d u^2)^2\right)$,
where the parameter $l$ is a real number. For $l=0$  the metric is Euclidean, while for $l=2$ it is the metric of the standard hyperbolic plane.
Those are the only two  manifolds of this family  which have constant curvature.
In \cite{Nishitani-Yagdjian} is  proved that the only stationary solutions of the equation (\ref{generalfns}) are the constant solutions and that
the global in time solvability can be destroyed
by parametric resonance phenomena. (For the scalar quasilinear wave equation
it was proved in  \cite{YagJMAA2001}.)
 For the parametric resonance phenomena in the scalar wave map-type  hyperbolic equations see \cite{yagdjian_birk} and references therein.
Then, according to  \cite{Ueda} (see  also   references therein)   the parametric resonance phenomena in the linear scalar wave equations can be localized in the space.
  \smallskip

Nakanishi and Ohta \cite{Nakanishi-Ohta} studied the Cauchy problem for the nonlinear wave equation
\begin{equation} \label{2}
\begin{cases}
    \square u  +  f(u) \left(\dot{u}^2  - |\nabla u|^2 \right)=0\,, \quad (t,x) \in {\mathbb R}^{1+n}\,, \cr
    u(0,x)= u_0 (x), \quad \dot{u}(0,x)= u_1 (x)\,, \quad x \in {\mathbb R}^{n}\,,
\end{cases}
\end{equation}
where $u=u(t,x) $ is a scalar real-valued unknown function, $f$ is a real valued smooth function.
In \cite{Nakanishi-Ohta} the following condition
\begin{equation}
        \label{1.3-Nakanishi-Ohta_K}
        \int_0^\infty \exp \left(\int_0^s f(r) dr \right) ds = \infty \quad \mbox{\rm and}  \quad  \int^0_{-\infty} \exp \left(\int_0^s f(r) dr \right) ds = \infty  
    \end{equation}
is suggested that is necessary and sufficient condition (Theorem 2.1~\cite{Nakanishi-Ohta}) for the existence of a global classical solution $u \in C^\infty ({\mathbb R}^{1+n}) $
for the  problem (\ref{2})  for any $u_0, u_1 \in C^\infty ({\mathbb R}^{n}) $.
Note here, that  the initial data $u_0, u_1$ are not assumed to be small.
The equation of (\ref{2}) is a model and special case  for wave maps.
\smallskip

The case of nonflat base manifold $L$ the wave maps are less investigated although they are of considerable interest in the   general relativity context.  
The Cauchy problem for  the wave maps in the perturbed Minkowski spacetime is considered in \cite{Choquet-Bruhat_ND} and \cite{Nishitani-Yagdjian} (cyclic universe).
More precisely, assume that $V= S\times {\mathbb R}$, with $S$ an $n$-dimensional orientable smooth manifold, and let $g$ be a Robertson-Walker metric
$
g   = -dt^2+a^2(t) \sigma$,
with the scale function   $a=a(t)$, where $\sigma =\sigma _{ij}\,dx^i \,dx^j $ is given smooth time independent metric on $S$, with non-zero injectivity radius.
\smallskip

Let $(S\times {\mathbb R},g)$ be a Robertson-Walker expanding universe with the metric
$
g   = -dt^2+a^2(t) \sigma,
$
 while $(S,\sigma ) $ is a smooth Riemannian manifold of dimension $n \leq 3$ with  non-zero injectivity radius   and $a=a(t)$ a positive increasing function of $t$ such that
$1/a(t)$ is integrable on $[t_0,\infty)$. Hence   a  domain of influence is permanently restricted (see, also, \cite[Sec.8]{yagdjian_birk}). 
Let $(M,h)$ be a proper Riemannian manifold regularly embedded in ${\mathbb R}^N$ such that \, Riem($h$) is uniformly bounded. 
Then according to Choquet-Bruhat~\cite{Choquet-Bruhat_ND} there exists a global wave map from $(S\times [t_0,\infty),g ) $ into  $(M,h)$ taking Cauchy data $\varphi  $, $ \psi $ with $D\varphi  $ and $\psi  $ in $H^1$
 if the integral of $1/a(t)$ on $[t_0,\infty)$ is less than some corresponding number $M(a,b)$.  
The number $M(a,b) $ depends on the initial data. Thus, (see Corollary on page 45~\cite{Choquet-Bruhat_ND}) under hypothesis of the theorem, for any finite value of the integral of $1/a(t)$ on $[t_0,\infty)$
there is an open set $U$ of   initial data in $H^1\times H^1$  such that if $(D \varphi , \psi ) \in U$, then there exists a global wave map taking the Cauchy data $(\varphi , \psi )$.
In particular, this is true for  the curved spacetime of the de~Sitter model of universe with the scale function   $a(t)=\exp(\Lambda t)$, $\Lambda >0 $.
\smallskip

D'Ancona and Zhang   \cite{D'Ancona}     derived the global existence of equivariant wave maps from the  so-called admissible manifolds to general targets for the small initial data of critical regularity. Both base  and target manifolds are assumed rotationally symmetric manifolds with global metrics 
\[
L\,:\, dr^2+g(r)^2d\omega ^2_{{\mathbb S}^{n-1}}\,,\qquad M\,:\, d\phi ^2+h(\phi )^2d\phi  ^2_{{\mathbb S}^{\ell-1}}\,,
\]
where $d\omega ^2_{{\mathbb S}^{n-1}}$ and $d\phi  ^2_{{\mathbb S}^{\ell-1}} $ are the standard metrics on the unit sphere. The solution has a form $u=( \phi ,\chi ) $ in coordinates on $M$, the 
radial component $\phi =\phi (t,r) $  depends only on time $t$ and $r$, the radial coordinate on $L$, while the
angular component $\chi  =\chi  (\omega )$ depends only on the angular coordinate $\omega  $ on $L$. Thus, $\chi \,:\,{\mathbb S}^{n-1}\longrightarrow   {\mathbb S}^{\ell-1}$ 
is a harmonic polynomial map of degree $k$, whose energy
density is $k(k+ n- 2)$ for some integer $k\geq  1$, while $\phi  $  satisfies the
{\it $\bar{\ell}$-equivariant wave map equation}
\begin{equation}
\label{1.7A}
\phi_{t t}-\phi_{r r}-(n-1) \frac{h^{\prime}(r)}{h(r)} \phi_{r}+\frac{\bar{\ell}}{h(r)^{2}} g(\phi) g^{\prime}(\phi)=0\,,
\end{equation}
where $\bar{\ell}=k(k+ n- 2) $. For (\ref{1.7A}) the authors  consider the Cauchy problem with initial data
$$
\phi(0, r)=\phi_{0}(r), \quad \phi_{t}(0, r)=\phi_{1}(r)\,.
$$
When $g(r)=r$ 
the problem for (\ref{1.7A}) reduces to the equation originally
studied in \cite{Shatah-Zadeh},\cite{Shatah-Struwe}. It is proved in \cite{D'Ancona} that on the admissible manifolds the wave flow satisfies smoothing and Strichartz estimates.
The metric $h$ of the base manifold is assumed to have a limit $  h^{\frac{1-n}{2}}(h^{\frac{n-1}{2}})^{\prime \prime }$ as $r \to \infty $.     
   The existence of small equivariant wave maps on admissible manifolds is proved in the critical space  
$H^{\frac{n}{2}}\times H^{\frac{n}{2}-1}$, and, moreover, the solution enjoys additional  $L^pL^q$
integrability properties  determined by the Strichartz estimates. 
  \smallskip

In the present paper we consider the wave map from the perturbed Minkowski spacetime, with the periodic in time perturbation,  into Riemannian manifold that is responsive (see self coherence structure below) to the parametric resonance generated by the metric $h$.
The result of the present note
requires some   assumption   on the ordinary
differential equation  related to the  parametric resonance generated by the periodic metric in $L$.
Consider the ordinary differential equation
\begin{eqnarray}
\label{Eastham1.2.1_new}
y_{tt}(t) + \left( \lambda    b^2(t)-q(t) \right) y(t)=0
\end{eqnarray}
with the periodic positive smooth non-constant function $b=b (t)$ and parameter $\lambda \in {\mathbb R} $. Let 
\begin{eqnarray*}
 q(t)
& =  &
 \frac{n}{4}\left(  \frac{n }{4}
-1 \right) \left( \frac{\dot{b}  (t)}{b(t)}\right)^2 - \frac{n}{2}  \frac{\ddot{b}   (t)}{b (t)}
\,.
\end{eqnarray*}
\noindent
{\bf Assumption ISIN }(\cite{Nishitani-Yagdjian}): {\it There exists the nonempty   open instability interval $\Lambda \subset (0,\infty)$ for equation (\ref{Eastham1.2.1_new}).} 
\medskip

We consider a wave map such that in the global chart of $M$ it can be written as a system of equations 
\begin{equation}
\label{4}
   u^i_{tt} -n \frac{\dot{b}(t)}{b(t)} u^i_t -  b^2(t) \Delta u^i
   + \sum_{j,k}\Gamma_{j,k}^i (u^1,\ldots,u^m) \left( {u_t}^j {u_t}^k- b^2(t) \nabla u^j \cdot \nabla u^k \right)=0,
\end{equation}
 $ i=1,\ldots, m$,  where $ b=b(t)$ is a smooth positive periodic function.
We are concerned with the small data global in time solution  to the Cauchy problem for equation (\ref{4}). 
Our main result shows that the global solvability is not a stable property under small perturbations of the wave map
if the Riemannian manifold $M$   
 {\it possesses a distinguished geodesic} (or {\it intrinsic self coherence structure}) in the sense of the following definition.

\begin{definition} \mbox{\rm \cite{yagdjian_geo}}
Riemannian or Lorentzian manifold $M$   
 {\it possesses a distinguished geodesic} (or {\it intrinsic self coherence structure}) if in some chart the straight half-line ${ \mathbb L}_+= \{(a_1t,\ldots,a_mt)\,|\, t \in (0,\infty) \}$ is covered by the geodesics.
\end{definition}
  The intrinsic self coherence structure can be characterized explicitly in the terms of Christoffel symbols $\Gamma_{j,k}^i $ as follows.

\begin{lemma}
\label{L1} \mbox{\rm \cite{yagdjian_geo}}
If in some chart of the Riemannian manifold  $M$ the segment $I $ of  the straight line  ${\mathbb L}= \{(a_1t,\ldots,a_mt)\,|\, t \in {\mathbb R} \}$ is covered by a smooth non-constant geodesic, then there is a function $f(t)$  such that 
\begin{equation}
\label{DL}
\sum_{j,k=1}^m \Gamma_{j,k}^i  (a_1t,\ldots,a_mt) a_ja_k = a_i f(t)  \,\, \, \mbox{\rm for all}\,\,\, t \in (a,b)  \subseteq {\mathbb R} \,\,\,\mbox{\rm and}\quad i=1,\ldots,m.
\end{equation}
Conversely, if  in some chart  there exists a continuously differentiable  function $f=f(t)$ such that (\ref{DL}) holds for all points of the segment $ I \subseteq {\mathbb L}$, then there is a geodesic covering    the segment $I $.
\end{lemma}
\medskip

The main result of this paper is  given by the following theorem.
\begin{theorem}
\label{MTH} Let  $b= b(t)$\,  be a defined on \,
${\R}$,\,  a  periodic, non-constant, smooth, and positive
function satisfying condition ISIT. Assume that the Riemannian manifold $M$ possesses  {\it intrinsic self coherence  structure} and for the function 
$ f(t)$,   $t \in {\mathbb R}$, 
the Nakanishi-Ohta condition (\ref{1.3-Nakanishi-Ohta_K}) does not hold, that is,
\begin{equation}
        \label{Not_N-O}
        \int_0^\infty \exp \left(\int_0^s f(r) dr \right) ds < \infty \quad \mbox{\rm or}  \quad  \int^0_{-\infty} \exp \left(\int_0^s f(r) dr \right) ds < \infty\,.
    \end{equation}  

Then for every
$n$, $s$, and  for every positive $\de$ there are initial data
$u_0^i , u_1^i \in C_0^\infty({\mathbb R}^n)$, $i=1,\ldots,m$, such that
\begin{equation}
\label{de}
\sum_{i=1}^m \|u_0^i \|_{(s+1)} + \|u_1^i \|_{(s)} \le \de\,,
\end{equation}
but the  solution
$u \in  C^2({\mathbb R}_+\times {\mathbb R}^n)$ to the problem with the prescribed data
\begin{equation}
\label{Cdata}
u^i(0,x) = u_0^i(x), \quad u_t^i(0,x) = u_1^i (x)\,, \quad i=1,\ldots, m, \quad x \in {\mathbb R}^n\,,
\end{equation}
for  the wave map   (\ref{4}) does not exist.
\end{theorem}

\begin{remark}
Assume that  $(u(s),0,\ldots,0)$ is geodesic   and for the function
\[
 f(t)=  \Gamma_{1,1}^1 (t,0,0,\ldots,0  ),     \quad   t \in {\mathbb R}_+,
\]
the Nakanishi-Ohta condition (\ref{1.3-Nakanishi-Ohta_K}) is not fulfilled. Then the statement of the theorem holds.
That is true also for any other coordinate axis.
\end{remark}

\begin{remark}
If (\ref{Sideris}) is fulfilled, 
then the system (\ref{4}) obeys intrinsic self coherence structure and the Nakanishi-Ohta condition (\ref{1.3-Nakanishi-Ohta_K}) is  fulfilled.   According to \cite{Sideris} the large data global solution exists for wave map without periodic perturbation ($b(t)\equiv 0$). The small amplitude solutions are known to exist (see,e.g.,\cite{Hormander}).  According to Theorem~\ref{MTH} (see, also, \cite{yagdjian_birk}) the     periodic perturbation $b(t)$ destroys global in time solvability even for the arbitrarily small data. 
\end{remark}

Following arguments of the proof Theorem 2.1~\cite{Nakanishi-Ohta} one can verify the  assertion of the next remark for the case of flat manifold 
 although we do not know if there is small data global existence for the case of non-flat $M$. 
\begin{remark}
The Cauchy problem for the system
\[
 u^i_{tt} -n \frac{\dot{b}(t)}{b(t)} u^i_{t }-  b^2(t) \Delta  u^i
   + f^i (u^i) \left(\left( {u_t}^i\right)^2- b^2(t) |\nabla u^i|^2 \right)=0, \quad i=1,\ldots, m\,,
\]
with the conditions (\ref{Cdata}) has a global solution $(u^1(x,t),\ldots,u^m(x,t)) \in C^\infty  $ for every \\
$(u^1_\ell(x ),\ldots,u^m_\ell(x )) \in C^\infty  ({\mathbb R}^n)\times \ldots \times C^\infty  ({\mathbb R}^n)$, $\ell=0,1$, 
 if and only if the condition 
\[
        \int_0^\infty \!\!\exp \left(\int_0^s \!f^i(r) dr \right) ds < \infty \,\,\, \mbox{\rm or}  \,\,\, \int^0_{-\infty}\! \!\exp \left(\int_0^s \!f^i(r) dr \right) ds < \infty\,,\,\,\, i=1,\ldots, m\,.
\]
is fulfilled.
\end{remark}
The proof of the next theorem is given in Section~\ref{S2}.
\begin{theorem}
\label{T4}
Let  $b= b(t)$\,be a defined on \,
${\R}$, periodic,  smooth, and positive
function. Assume that the Riemannian manifold $M$ possesses    intrinsic self coherence  structure  and the Cauchy problem for (\ref{4}) has a global solution $(u^1(x,t),\ldots,u^m(x,t)) \in   C^2({\mathbb R}_+\times {\mathbb R}^n)$ for every  initial data 
$(u^1_i(x ),\ldots,u^m_i(x )) \in C^\infty  ({\mathbb R}^n)\times \ldots \times C^\infty  ({\mathbb R}^n)$, $i=0,1$. Then 
  the Nakanishi-Ohta condition (\ref{1.3-Nakanishi-Ohta_K}) is  fulfilled.    
\end{theorem}
Note   that  the initial data $u_0, u_1$ are not assumed small.
Existence of the distinguished geodesics allows also to extend result of \cite{Nakanishi-Ohta} from the wave map type  equations to the  wave map  with the non-oscillating coefficients for some non-small initial data. That will be proved in the forthcoming paper. 
\medskip

The present paper is organized as follows. In Sec.~\ref{S1} we illustrate   Theorem~\ref{MTH} by several examples. 
Then, in Sec.~\ref{S2}, we lower  the system of equations to the single scalar equation. In Sec.~\ref{S3} we 
describe some elements of Floquet-Lyapunov theory  with its application to the parametric resonance in the ordinary differential equations. In Sec.~\ref{S4} and Sec.~\ref{S5} 
we complete the proofs of   Theorem~\ref{MTH} and Theorem~\ref{T4}, respectively. The final Sec.~\ref{S6} is devoted to the proof of Lemma~\ref{L1}.

\section{Illustration of Theorem~\ref{MTH} by Examples}
\label{S1}

In the  spacetime with the metric tensor 
\begin{eqnarray*}
g_{ik}
& = &
\begin{pmatrix}
1 & 0 & 0& 0 \cr
0 & -a^2(t)& 0& 0 \cr
0 & 0 & -a^2(t)& 0 \cr
\vdots & \vdots & \ddots  & \vdots \cr
0 &  0 & 0&  -a^2(t)\cr 
\end{pmatrix}, \quad |g|= a^{2n}(t), \\
\end{eqnarray*}
the covariant D'Alembert operator  is defined as follows:  
\begin{eqnarray*}
\Box_g u= \frac{1}{\sqrt{|g|}}\frac{\partial}{\partial x^i}\left( \sqrt{|g|}g^{ik} \frac{\partial}{\partial x^k} u\right)= \frac{\partial^2}{\partial t^2}u
+ n\frac{\dot{a} (t)}{a(t)}\frac{\partial}{\partial t}u-\frac{1}{a^2(t)} \Delta u \,.
\end{eqnarray*}
If we denote $b(t)=1/a(t)$, then
\begin{eqnarray*}
\Box_g u=   \frac{\partial^2}{\partial t^2}u
- n\frac{b'(t)}{b(t)}\frac{\partial}{\partial t}u- b(t)^2  \Delta u\,.
\end{eqnarray*}
The corresponding wave map equation is (\ref{4}). 
Cyclic spacetime with the periodic smooth positive scale factor $a=a(t)$ is one of the models  of the cosmology (see \cite[Ch. 9]{Ohanian-R}).
\medskip

\noindent
\textbf{Example 1:} Consider the system (\ref{4}) with $m=2$:
\begin{eqnarray}
\label{14}
  \begin{cases} 
  \left(  \partial_t^2 -n \frac{\dot{b}(t)}{b(t)}\partial_t-  b^2(t) \Delta \right)u^1 \cr
\hspace{2cm} + \sum_{j,k=1}^2\Gamma_{j,k}^1 (u^1,u^2) \left(\dot{u}^j\dot{u}^k- b^2(t) \nabla u^j \cdot \nabla u^k \right) =0,  \cr
    \left(  \partial_t^2 -n \frac{\dot{b}(t)}{b(t)}\partial_t-  b^2(t) \Delta \right)u^2 \cr
\hspace{2cm}+ \sum_{j,k=1}^2\Gamma_{j,k}^2 (u^1,u^2) \left(\dot{u}^j\dot{u}^k- b^2(t) \nabla u^j \cdot \nabla u^k \right) =0.
    \end{cases}
\end{eqnarray}
We define in $M $  the diagonal metric tensor $h_{ik} (u^1,u^2): = h(u^1,u^2)\delta_{ik}$.
Then, the Christoffel symbols are:
\[
    \Gamma_{j,k}^i =
\frac{1}{2h(u^1,u^2)} \left(\frac{\pa}{\pa u^j}h(u^1,u^2)\delta_{ki} + \frac{\pa}{\pa u^k}h(u^1,u^2)\delta_{ji} - \frac{\pa}{\pa u^i}h(u^1,u^2)\delta_{kj} \right),
\]
where $i,j,k=1,2 $. Hence,
\begin{eqnarray*}
 &  &
   \Gamma_{1,1}^1 =-\Gamma_{2,2}^1   = \Gamma_{2,1}^2= \Gamma_{1,2}^2 = \frac{ 1}{2h(u^1,u^2)} \left( \frac{\pa}{\pa u^1}h(u^1,u^2) \right), \\
  &  &
  \Gamma_{2,1}^1=\Gamma_{1,2}^1  = \Gamma_{2,2}^2 =- \Gamma_{1,1}^2=\frac{ 1}{2h(u^1,u^2)} \left( \frac{\pa}{\pa u^2}h(u^1,u^2) \right).
\end{eqnarray*}
The Gaussian curvature of the surface with such metric is
\[
K=-\frac{1}{h(u^1,u^2)}\Delta \ln h(u^1,u^2)\,.
\]
The wave  map equation (\ref{4}) reads 
\begin{eqnarray*}
&  &
 \begin{cases}
   \left( \dsp  \partial_t^2 -n \frac{\dot{b}(t)}{b(t)}\partial_t-  b^2(t) \Delta \right)u^1\cr
\hspace{2cm} +
\dsp \frac{ 1}{2h(u^1,u^2)} \left( \frac{\pa}{\pa u^1}h(u^1,u^2) \right)    (\dot{u}^1\dot{u}^1- b^2(t) \nabla u^1 \cdot \nabla u^1) \cr
\hspace{2cm} \dsp + \frac{ 1}{ h(u^1,u^2)} \left( \frac{\pa}{\pa u^2}h(u^1,u^2) \right)  (\dot{u}^1\dot{u}^2- b^2(t) \nabla u^1 \cdot \nabla u^2) \cr
\hspace{2cm} \dsp -\frac{ 1}{2h(u^1,u^2)} \left( \frac{\pa}{\pa u^1}h(u^1,u^2) \right)   (\dot{u}^2\dot{u}^2- b^2(t) \nabla u^2 \cdot \nabla u^2) =0,  \cr
    \left(  \dsp \partial_t^2 -n \frac{\dot{b}(t)}{b(t)}\partial_t-  b^2(t) \Delta \right)u^2 \cr
\hspace{2cm}-\dsp \frac{ 1}{2h(u^1,u^2)} \left( \frac{\pa}{\pa u^2}h(u^1,u^2)  \right)  (\dot{u}^1\dot{u}^1- b^2(t) \nabla u^1 \cdot \nabla u^1) \cr
\hspace{2cm} \dsp + \frac{ 1}{ h(u^1,u^2)} \left( \frac{\pa}{\pa u^1}h(u^1,u^2) \right)  (\dot{u}^1\dot{u}^2- b^2(t) \nabla u^1 \cdot \nabla u^2) \cr
\hspace{2cm} +\dsp \frac{ 1}{2h(u^1,u^2)} \left( \frac{\pa}{\pa u^2}h(u^1,u^2) \right)   (\dot{u}^2\dot{u}^2- b^2(t) \nabla u^2 \cdot \nabla u^2) =0\, . 
    \end{cases}
\end{eqnarray*}
If $b(t)=const>0$, the small amplitude solutions of (\ref{14}) exist globally. 
Now we focus on the case with a half-diagonal  ${\mathbb L}_+=  \{( t,\ldots, t)\,|\, t \in (0,\infty) \}\subset {\mathbb D}$. We note that 
\[
    \sum_{j,k=1}^2 \Gamma_{jk}^1(u^1,u^2)
         = \frac{ 1}{h} \left( \frac{\pa}{\pa u^2}h\right),\qquad
    \sum_{j,k=1}^2 \Gamma_{jk}^2(u^1,u^2)
        = \frac{ 1}{h} \left( \frac{\pa}{\pa u^1}h\right).
\]
Assume that
\[
\frac{\pa h}{\pa u^k}(u^1,u^2) =  \frac{\partial h}{\partial u^l}(u^1,u^2)
\quad {\rm if} \quad u^1 = u^2  \quad {\rm for} \quad k,l = 1,2.
\]
Then, due to the last assumption  on $h_{ik}$ we  set $a_1= a_2=1 $ and obtain the function of (\ref{DL}) 
\[
    f(\xi) :=  
\sum_{j,k=1}^2 \Gamma_{jk}^1 (\xi, \xi)= \sum_{j,k=1}^2 \Gamma_{jk}^2(\xi, \xi) \quad
   \mbox{\rm if}\quad   \xi \in \mathbb{R}_+\,.
\]

To find geodesics
let $(U, \varphi)$ be a  parametrization of the  manifold  $M$ and let $\alpha : I \rightarrow M$
be a curve parametrized by arc length, whose trace is contained in $\varphi(U) .$ Write
\[
\alpha(s)=\varphi(u(s), v(s))\,,
\]
where $ u= u(s)$ and $v= v(s)$ are real-valued functions of $ s $.   Then $\alpha$ is a geodesic if 
\[  
\begin{cases} 
\dsp \ddot{u}(s)
  +  
\frac{ 1}{2h(u ,v)} \left( \frac{\pa}{\pa u}h(u,v) \right)\left(\dot{u}(s)\right)^{2}\cr
\dsp \hspace{1.4cm} +  \frac{ 1}{ h(u ,v)} \left( \frac{\pa}{\pa v}h(u ,v) \right) \dot{u}(s)  \dot{v}  (s)
-\frac{ 1}{2h(u ,v)} \left( \frac{\pa}{\pa u}h(u,v) \right)\left(\dot{v} (s)\right)^{2}=0\,,   \cr
\dsp \ddot{v}(s)
  -  
\frac{ 1}{2h(u ,v)} \left( \frac{\pa}{\pa v}h(u ,v) \right)\left(\dot{u}(s)\right)^{2}\cr
\dsp \hspace{1.4cm} +  \frac{ 1}{ h(u ,v)} \left( \frac{\pa}{\pa u }h(u ,v) \right) \dot{u}(s) \dot{v} (s)+\frac{ 1}{2h(u ,v)} \left( \frac{\pa}{\pa v}h(u ,v) \right)\left(\dot{v} (s)\right)^{2}=0\,. 
\end{cases}
\]  
We claim that there exists a geodesic curve that lies in the diagonal  $\mathbb D $. Indeed, set $u(s)=v(s) $, then equation of geodesic and unite speed equation read
\begin{eqnarray*}
& &
\ddot{u} (s)
+  \frac{ 1}{ h(u(s) ,u(s))} \left( \frac{\pa}{\pa u}h(u(s) ,u(s)) \right) (\dot{u} (s))^2=0 \,,  \\
& &
1  =
h(u(s) ,u(s))2(\dot{u} (s))^2\,.
\end{eqnarray*}
From the second equation 
the solution $u=u(s)$ can be given implicitly by
\begin{eqnarray}
\label{geo}
&  &
\int_0^{u(s)}\sqrt{h(r  ,r )}  \,dr= \frac{1}{\sqrt{2}} s+C\,.
\end{eqnarray}
Let $ h(u^1,u^2)=(1+u_1^2+u_2^2)^\alpha $, we check condition  (\ref{1.3-Nakanishi-Ohta_K}):
\begin{eqnarray*}
&  &
\int_0^{\pm \infty} \exp \left(\int_0^s f(r) dr \right) ds
= 
 \int_0^{\pm \infty} ( 1+2s^2  )^\alpha  ds
= \int_0^{\pm \infty}  h(s,s)  ds
={\pm \infty}  \,.
\end{eqnarray*}
Hence, the condition (\ref{1.3-Nakanishi-Ohta_K}) is equivalent to  the inequality $\alpha >-\frac{1}{2}  $.
For the case of $ h(u ,v)=(1+u ^2+v^2)^\alpha $ the equation (\ref{geo}) for the geodesics   leads to 
the function  $u=u(s)$ that is defined  implicitly by 
\begin{equation}
\label{15new}
u  F \left(\frac{1}{2},-\frac{\alpha }{2} ;\frac{3}{2};-2 u^2\right)= \frac{1}{\sqrt{2}}s+C\,.
\end{equation}
If $\alpha =-1 $, then the condition (\ref{1.3-Nakanishi-Ohta_K}) is violated and the equation (\ref{15new})  simplifies to $ u(s)= C_1e^s+ C_2e^{-s}$ that implies for the geodesic 
\[
  u(s)=v(s)=  C_1e^s+ C_2e^{-s}\,.
\]
The non-constant geodesic that    belongs to the diagonal $\mathbb D $ and starts at the origin is given by 
\[
u(s)=v(s)=
\frac{1}{ \sqrt{2}}\sinh (s)\,.
\]
For the case of $ h(u^1,u^2)=(1+u_1^2+u_2^2)^{-1} $ on the diagonal $\mathbb D$ the Christoffel symbols are  
\[
   \Gamma_{1,1}^1 
 =
-\Gamma_{2,2}^1   = \Gamma_{2,1}^2= \Gamma_{1,2}^2 = 
  \Gamma_{2,1}^1  
  =  
\Gamma_{1,2}^1  
   =   \Gamma_{2,2}^2 =- \Gamma_{1,1}^2
=  -\frac{1}{ \sqrt{2}}\tanh (s) \sech (s)  . 
\]
The Gaussian curvature of the surface with the metric $ h(u^1,u^2)=(1+u_1^2+u_2^2)^\alpha $ is
\begin{eqnarray*}
K
& = &
-\frac{1}{h(u^1,u^2)}\Delta \ln h(u^1,u^2)=
-4\alpha(1+u ^2+v^2)^{-\alpha-2}\,.
\end{eqnarray*}
It is also a scalar curvature. It is constant iff $\alpha =-2 $.

\smallskip

\noindent
\textbf{Example 2:} Define the metric  $ h(u ,v)=  (1+v)^{-\ell}$, $\ell \geq 0 $ on  $M=\{(u ,v)\in {\mathbb R}\,|\, v>-1\} $,  then the Christoffel symbols are 
\begin{eqnarray*} 
  \Gamma^1_{2, 1} = \Gamma^1_{1,2 } = \Gamma^2_{2, 2}=-\Gamma^2_{1, 1}  &=&  -\frac{\ell }{2 (1+v)}
\end{eqnarray*}
while the equations for the geodesics are
\begin{eqnarray*}
\begin{cases} 
\dsp \ddot{u}(s)     -\frac{\ell }{ (1+v)}   \dot{u}(s)  \dot{v}  (s) =0  \,,\cr
\dsp \ddot{v}(s)
+\frac{\ell }{2 (1+v)}\left(\dot{u}(s)\right)^{2}   -\frac{\ell }{2 (1+v)}\left(\dot{v} (s)\right)^{2}=0\,. 
\end{cases}
\end{eqnarray*} 
If $\ell=2 $ this system has a solution $u(s)=u(0)$, $v(s)=Ce^s-1$, that is a vertical half-line in the positive half-plane. The geodesic starting at the origin is $u(s)=0 $, $v(s)=e^s-1$. Then, 
\[
f(t)=  -\frac{\ell }{2 (1+t)}, \quad   \int_0^\infty \exp \left(\int_0^s f(r) dr \right) ds = 
\int_0^\infty  (1+s)^{-\frac{\ell}{2}} ds< \infty
\]
implies $\ell>2  $. For the case of $\ell \in[0,2) $ the nonexistence of the global solution for arbitrary small data is proved in \cite{Nishitani-Yagdjian}. The  global existence  of arbitrary small data solutions for the case of  $\ell=2 $ and non-constant periodic $b=b(t) $ remains an open problem.

\smallskip

\noindent
\textbf{Example 3:} Assume now that $ h(u^1,u^2)=(1+u_1^2+u_2^4)^\alpha =(1+u ^2+v^4)^\alpha $, then the Christoffel symbols are 
\begin{eqnarray*} 
  \Gamma^1_{1, 1}=-\Gamma^1_{2, 2} = \Gamma^2_{2, 1} = \Gamma^2_{1,2 }  &=& \frac{\alpha  u}{u^2+v^4+1} \,,\\
 \Gamma^1_{2, 1} = \Gamma^1_{1,2 } = \Gamma^2_{2, 2}=-\Gamma^2_{1, 1}   &=& \frac{2 \alpha  v^3}{u^2+v^4+1}\,, 
\end{eqnarray*} 
and the equations for the geodesics are
\begin{eqnarray*}
\begin{cases} 
\dsp \ddot{u}(s)
  +  
\frac{  \alpha u}{ (1+u ^2+v^4) } \left(\dot{u}(s)\right)^{2}\cr
\dsp \hspace{2cm} +  \frac{ 4v^3 \alpha }{ (1+u ^2+v^4) }   \dot{u}(s)  \dot{v}  (s)
-\frac{  \alpha u}{ (1+u ^2+v^4) }  \left(\dot{v} (s)\right)^{2}=0   \,,\cr
\dsp \ddot{v}(s)
  -  
 \frac{ 2v^3 \alpha }{ (1+u ^2+v^4) }\left(\dot{u}(s)\right)^{2}\cr
\dsp \hspace{2cm} +  \frac{ 2\alpha u }{ (1+u ^2+v^4) }  \dot{u}(s) \dot{v} (s)+\frac{ 2v^3 \alpha }{ (1+u ^2+v^4) }\left(\dot{v} (s)\right)^{2}=0\,. 
\end{cases}
\end{eqnarray*}
The curve $v (s)=0$ is geodesic if 
\begin{eqnarray*} 
\dsp \ddot{u}(s)
  +  
\frac{  \alpha u(s)}{ (1+u ^2(s)) } \left(\dot{u}(s)\right)^{2} =0, \quad 1  =
h(u(s) ,u(s))(\dot{u} (s))^2  \,,
\end{eqnarray*}
that is,  
\begin{eqnarray*} 
\dsp \ddot{u}(s)
  +  
\frac{  \alpha u(s)}{ (1+u ^2(s)) } \left(\dot{u}(s)\right)^{2} =0, \quad 1  =
(1+u ^2(s))^\alpha(\dot{u} (s))^2  \,.
\end{eqnarray*}
The function $f(t)=  \alpha t/(1+t^2) $ and
\[    
\int_0^\infty\! \exp \left(\int_0^s f(r) dr \right) ds = 
\int_0^\infty \! \exp\left(\int_0^s \frac{\alpha r}{1+r^2} dr \right) ds 
=  \int_0^\infty  (1+s^2)^{\alpha/2 } ds< \infty\,.
\]
The condition (\ref{1.3-Nakanishi-Ohta_K}) implies $\alpha >-1  $. 

The line $u(s)=0$ is also a geodesic and with the function $f(t)= 2\alpha t^3/(1+t^4) $ the condition (\ref{1.3-Nakanishi-Ohta_K}) 
\[    
\int_0^\infty \exp \left(\int_0^s f(r) dr \right) ds = 
\int_0^\infty \exp \left(\int_0^s \frac{2\alpha r^3}{1+r^4} dr \right) ds 
=  \int_0^\infty  (1+s^4)^{\alpha/2 } ds< \infty
\]
reads $\alpha >-1/2 $. Thus, the choice of the geodesic line is essential. 
The Gaussian curvature of the surface with the metric $ h(u^1,u^2)=(1+u_1^2+u_2^4)^\alpha =(1+u ^2+v^4)^\alpha $ is
\begin{eqnarray*}
K
& = &
-2 \alpha  \left(u^2 \left(6 v^2-1\right)-2 v^6+v^4+6 v^2+1\right) \left(u^2+v^4+1\right)^{-\alpha -2}\,.
\end{eqnarray*} 
It is also a scalar curvature. 
\smallskip

The next example shows that small perturbation of the diagonal metric tensor does not  eliminate blow up phenomenon.  
\smallskip

\noindent
\textbf{Example 4:} Let ${\mathbb R}^m $ be provided with the metric defined by the metric tensor $h_{ik} (u ) = h (u )(\delta_{ik} + H_{ik}(u ))$,   where $u =(u^1, \dots , u^m)$ and $h=h (u ) $ is smooth  positive  function. We denote  $M$ such   Riemannian manifold. Assume that $H(u)$ is a smooth matrix function  with the matrix norm $\|H(u)\|<1$     and that on the diagonal ${\mathbb D}$ of $M$ 
\begin{eqnarray*}
&  &
\frac{\pa}{\pa u^k} H( u)=0 , \quad H( u)=0  \quad  { \rm if } \quad  u \in {\mathbb D},\quad  \forall k =1,2,\ldots,m,\\
&  &
    \frac{\pa}{\pa u^k} h(u^1, \dots , u^m) = \frac{\pa}{\pa u^l} h(u^1, \dots , u^m) 
 \quad {\rm if } \quad  u \in {\mathbb D},\quad  \forall k,l=1,2,\ldots,m.
\end{eqnarray*}
The Christoffel symbols for the metric $h_{ik} (u)$ on the diagonal $\mathbb D$  are:
\begin{eqnarray*}
 \Gamma^i_{jk} (u)
&= &
\frac{1}{2}    \frac{1}{h(u)}   \left( \frac{\pa}{\pa u^j}  h(u) \delta_{k i}  + \frac{\pa}{\pa u^k} h (u)\delta_{j i}  -  \frac{\pa}{\pa u^i}  h(u) \delta_{jk}  \right)
\end{eqnarray*}
and
\begin{eqnarray*}
\sum_{j,k=1}^m \Gamma^i_{jk} (u)
&= &
\frac{1}{2} m   \frac{1}{h(u)}  \frac{\pa}{\pa u^1}  h(u) ,\quad i=1,\ldots,m,\quad    u\in {\mathbb D}\,.
\end{eqnarray*}
The diagonal $\mathbb D $ is a geodesic. Indeed, we set the initial conditions 
\[
u^1(0)=\ldots=u^m(0)=0\,,\qquad  \frac{d  u^1}{ds } (0)=\ldots=\frac{d  u^m}{ds }(0)=\left( m h  ( 1   ,\ldots,1   ) \right)^{-1/2} \,,
\]
and consider the function $\tilde{u} =\tilde{u}(s)  $ that  solves the Cauchy problem 
\[
\frac{d^2 \tilde{u} }{ds ^2}  
+ \frac{1}{2} m   \frac{1}{h(u)}  \frac{\pa}{\pa u^1}  h(u)  \left( \frac{d  \tilde{u} }{ds } \right)^2=0\,,  
\qquad  \tilde{u}(0)=0,\quad \frac{d   \tilde{u} }{ds }(0)=\left( m h  ( 1   ,\ldots,1   ) \right)^{-1/2} . 
\]
Then the function $u(s)= (\tilde{u}(s),\ldots,\tilde{u}(s))$ is   a geodesics that lies in $\mathbb D $. 
Therefore, if we define
\[
    f(u) := \frac{m}{2h(u)}     \frac{\pa}{\pa u^1}  h(u)  ,\quad    u\in \mathbb{D},
\]
then with $a_1=\ldots=a_m=1 $ the condition (\ref{DL}) is fulfilled:
\[
    \sum_{j,k=1}^m \Gamma^1_{jk}(u) = \sum_{j,k=1}^m \Gamma^2_{jk}(u) = \ldots = \sum_{j,k=1}^m \Gamma^m_{jk}(u) = f(u),\quad    u\in \mathbb{D}.
\]
In order to verify the condition (\ref{Not_N-O}) we   specify $h(u)=(1+u_1^2+ \ldots  +u_m^2)^\alpha  $, then
\[
    f(u) := \frac{m\alpha u}{1+mu ^2}    ,\quad    u\in {\mathbb R},
\]
\[    
\int_0^\infty \! \!\exp \left(\int_0^s \!\!f(r) dr \right) ds = 
\int_0^\infty\!\! \exp \left(\int_0^s \!\! \frac{m\alpha r}{1+mr^2} dr \right) ds 
=  \int_0^\infty \!\! (1+ms^2)^{\alpha /2} ds< \infty\,.
\] 
Condition (\ref{Not_N-O}) implies $\alpha < -1  $. 

\medskip

\noindent
\textbf{Example 5:} Let  $b(t)  =\sqrt{ 1 + \varepsilon \sin{(t)}}$,\, where $\varepsilon \in(0,1) $,   be a defined on \, ${\R}$,\,  a  periodic, non-constant, smooth, and positive
function. Assume that $m=2$, then 
\[
 \begin{cases} 
\displaystyle  \left(  \partial_t^2 -n \frac{\varepsilon \cos(t)}{2(1+ \varepsilon \sin(t) )}\partial_t-    (1+ \varepsilon \sin(t) )\Delta \right)u
   +    |\dot{v}|^2 - (1+ \varepsilon \sin(t) )|\nabla v|^2    =0,\cr
\displaystyle    \left(  \partial_t^2 -n \frac{\varepsilon \cos(t)}{2(1+ \varepsilon \sin(t) )}\partial_t-    (1+ \varepsilon \sin(t) )\Delta \right)v
   +    |\dot{u}|^2 - (1+ \varepsilon \sin(t) )|\nabla u|^2    =0.  
   \end{cases}
\]
Then for every $n$, $s$, and  for every positive $\de$ there are data $u_0 , v_0,  u_1 ,v_1 \in C_0^\infty({\R}^n)$ such that
\begin{eqnarray*}
  \|u_0  \|_{(s+1)} + \|u_1  \|_{(s)}+ \|v_0  \|_{(s+1)} + \|v_1  \|_{(s)} \le \de
\end{eqnarray*}
but the solution $u, v \in  C^2({\R}_+\times {\R}^n)$ to the problem with data
\[
u (0,x) = u_0 (x),  \,\,\,  u_t (0,x) = u_1  (x),  \,\,\, v (0,x) = v_0 (x),   \,\,\,
v_t (0,x) = v_1  (x)\,,     \,\,\, x \in {\R}^n
\]
does not exist. For the same data if $\varepsilon  = 0$ then a small data solution exists globally.
The Riemannian curvature of this spacetime with $n=3$ is 
\[
-\frac{3 \varepsilon  (\varepsilon  \cos (2 t)+3\varepsilon +2 \sin (t))}{2 (\varepsilon  \sin (t)+1)^2}\,,
\] 
which is sign changing  in time.

\medskip

\section{Lowering to the scalar equation}
\label{S2}

 The main idea is to use a composition of the solution of the wave equation in $L$ with the distinguished geodesic of the target manifold  $M$.    This composition is a wave map. For the properly chosen geodesic  such wave map   blows up for the large time (see also \cite{Nishitani-Yagdjian}).  
 Consider the system of equations
\begin{eqnarray*}
   u_{tt}^i-n \frac{\dot{b}(t)}{b(t)} u_{t }^i-  b^2(t) \Delta  u^i
   + \sum_{j,k}\Gamma_{j,k}^i (u^1,\ldots,u^m) \left( {u}^j_t {u}^k_t- b^2(t) \nabla u^j \cdot \nabla u^k \right)=0,
\end{eqnarray*}
$i=1,\ldots m$, where $\Gamma_{j,k}^i (u) $, $b(t)$  are $C^\infty$ functions satisfying condition  (\ref{DL}). The choice
of the initial data 
\[
u^i(0,x) =a_iu_0(x), \quad u_t^i(0,x) = a_iu_1 (x)\,, \quad \quad i=1,\ldots m, \quad x
\in {\mathbb  R}^n\,,
\]
for the system of equations and the  intrinsic self coherent structure of the manifold force a unique local solution to be on the track of the  distinguished geodesic.   This allows the lowering of the wave map system to the scalar equation. Indeed, if we  consider   the Cauchy problem for the auxiliary scalar   equation 
\begin{eqnarray}
\label{15}
\label{20}
\begin{cases}  
\dsp   u_{tt} -n \frac{\dot{b}(t)}{b(t)} u_{t }-  b^2(t) \Delta  u 
   +  f(u) \left(    u_t^2 - b^2(t) \nabla u  \cdot \nabla u  \right)=0, \cr
u (0,x) = u_0(x), \quad u_t (0,x) =  u_1 (x)\,, \quad   x
\in {\mathbb  R}^n\,, 
\end{cases} 
\end{eqnarray}
then  according to the uniqueness of the solution we have 
\[
u^1(t,x) =a^1u(t,x),\, u^2(t,x) =a^2u(t,x),\,\ldots ,\,u^m(t,x)=a^mu(t,x)    
\]
for all $x \in   {\mathbb  R}^n$, \,  $t \geq 0$. Thus  we can restrict ourselves  to   the Cauchy problem (\ref{15}) for the auxiliary scalar   equation,  
where $f (u) $, $b(t)$  are $C^\infty$ functions and  $f (u) $ is from condition  (\ref{DL}). 
For this Cauchy problem we find arbitrarily small smooth initial  data and prove that the solution blows up in finite time. This implies that the solution to the problem (\ref{4})\&(\ref{Cdata}) blows up in finite time that proves  Theorem~\ref{MTH}.

Consider the   equation of (\ref{15}). 
By the Hopf-Cole-Nakanishi-Ohta transformation 
\begin{equation}
\label{1.4-Nakanishi-Ohta-K}
v = G(u) :=  \int_0^u \exp \left( \int_0^s f(r) dr \right) ds \,,
\end{equation}
the equation   (\ref{20}) is transformed into the linear wave equation 
\begin{equation}
\label{lintrans}
    v_{tt} - n\frac{\dot{b}(t) }{b(t)}v_t- b^2(t)\De v = 0\,.
\end{equation}
Since $G \in C^2({\mathbb R})$ and $G'>0$, there exists the inverse of $G$:
\begin{equation}
\label{H}
H:= G^{-1} \in C^2(a,b),
\end{equation}
where we denote
\begin{equation}
\label{ab}
a:= \lim_{u \to -\infty} G(u)\,, \qquad b := \lim_{u \to \infty} G(u)\,.
\end{equation}

Next we apply the partial Liouville  transformation that eliminates the first derivative $v_t $ in 
(\ref{lintrans}). 
 More precisely, we set   
\begin{eqnarray*}
&  &
v= b^{\frac{n}{2}}(t)w\,,\quad b(t)=1/a(t), 
\end{eqnarray*}
then
\begin{eqnarray*}
&  &
 v_{tt}- n \frac{\dot{b}(t) }{b(t)} v_t- b^2(t) \Delta v\\
& = &
b^{\frac{n}{2}}(t)\Bigg[ w_{tt}- b^{2}(t) \Delta w  +\left\{ \frac{n}{2}\left(1-\frac{n}{2}  \right)  
 \left( \frac{d}{d t}\frac{1}{b(t)}\right)^2 b^2(t)- \frac{n}{2}   \left( \frac{d^2}{d t^2}\frac{1}{b(t)}\right) b(t)  \right\} w  \Bigg].
\end{eqnarray*}
Thus, we have to study the following linear hyperbolic equation 

\begin{eqnarray*}
 w_{tt}- b^{2}(t) \Delta w  +\left(  \frac{n}{2}\left(1-\frac{n}{2}  \right)  
 \left( \frac{d}{d t}\frac{1}{b(t)}\right)^2 b^2(t)- \frac{n}{2}   \left( \frac{d^2}{d t^2}\frac{1}{b(t)}\right) b(t) \right) w 
=0 
\end{eqnarray*}
with  the 1-periodic  positive smooth function $b=b(t)$.

\section{Floquet-Lyapunov theory. Parametric Resonance in ODE}
\label{S3}

We are going to apply the Floquet-Lyapunov theory for the ordinary differential equation with the periodic coefficients. 
Consider the ordinary differential equation:
\begin{eqnarray*}
 W_{tt}+\left( \lambda   b^{2}(t)  + \frac{n}{2}\left(1-\frac{n}{2}  \right)  
 \left( \frac{d}{d t}\frac{1}{b(t)}\right)^2 b^2(t)- \frac{n}{2}   \left( \frac{d^2}{d t^2}\frac{1}{b(t)}\right) b(t) \right) w 
=0 
\end{eqnarray*}
with the periodic positive smooth non-constant function $b=b(t)$ and parameter $\lambda \in {\mathbb R} $.

It is more convenient  to rewrite this equation by means of the new positive periodic function   
\[
\alpha (t):=  b^{2}(t),  
\]
then 
\begin{eqnarray*}
 W_{tt}+\left\{ \lambda   \alpha (t) - \frac{n}{4}\left[ \frac{3}{2} \left( \frac{\dot{\alpha }  (t)}{\alpha (t)}\right)^2- \frac{\ddot{\alpha }  (t)}{\alpha (t)} \right]  
- \frac{n}{8}\left(\frac{n}{2}- 1 \right) \left( \frac{\dot{\alpha }(t)}{\alpha (t)}\right)^2  \right\} W 
=0 \,.
\end{eqnarray*}
Consider now the equation
\begin{eqnarray}
\label{Eastham1.2.1}
y_{tt}(t) + \left( \lambda    \alpha (t)-q(t) \right) y(t)=0
\end{eqnarray}
with  the periodic coefficients $\alpha (t)
  =   
b^{ 2}(t) $ and  
\begin{eqnarray*} 
 q(t)
& =  & 
\frac{n}{4}\left[ \frac{3}{2} \left( \frac{\dot{\alpha }  (t)}{\alpha (t)}\right)^2- \frac{\ddot{\alpha }  (t)}{\alpha (t)} \right]  
- \frac{n}{8}\left(\frac{n}{2}- 1 \right) \left( \frac{\dot{\alpha }(t)}{\alpha (t)}\right)^2  \,.
\end{eqnarray*}
The first part of the last expression is the so-called Schwarz derivative for the antiderivative of $\alpha (t)$. For equation (\ref{Eastham1.2.1}) the spectrum of the  eigenvalue problem with the boundary condition 
\[
y(0)=y(1)=0
\]
is discrete. 
The equation (\ref{Eastham1.2.1})  can be written also as a system of differential equations for the
vector-valued function $x(t)={}^t(w_t,w)$:
\[
\displaystyle{\frac{d}{dt}} x(t) = A(t)
 x(t)\,,\qquad  \mbox{\rm where}
\quad A(t) := \left(\begin{array}{lll}
0  & - \lambda    \alpha (t)+q(t) \\
1  & 0 
\end{array} \right)\,.
\]
Let the  matrix-valued function\, $X_\lambda (t,t_0)$,\, depending on \,$\lambda $,\,
be a solution of the 
Cauchy problem
\begin{equation}
\label{system}
\displaystyle{\frac{d}{dt}} X = A(t)
 X\,,\qquad 
X(t_0,t_0) = 
\left(\begin{array}{lll}
1  & 0 \\
0  & 1 
\end{array} \right)\,.
\end{equation}
Thus, \, $X_\lambda (t,t_0)$\, gives a fundamental solution  to the 
equation (\ref{Eastham1.2.1}).  
In what follows we often omit subindex $\lambda $ of $X_\lambda  (t,t_0)$. 
The Liouville formula 
\begin{eqnarray*}
W(t) = W(t_0) \exp \left( \int_{t_0}^t S(\tau) d\tau \right),\\ 
\end{eqnarray*}
 where $
W(t):= \det X(t,t_0)$, \,
$S(t):= \sum_{k=1}^2 A_{kk}(t)$ 
with $S(t) \equiv 0$ guarantees the existence of the inverse matrix
$X_\lambda  (t,t_0)^{-1}$. 
 For the  matrix $ X(1,0) $ we
will use a notation
\[
X_\lambda (1,0) =  \left(\begin{array}{lll}
 b_{11}  & b_{12}   \\
b_{21}      & b_{22}  
\end{array} \right)\,.
\]
This matrix is called a {\it monodromy matrix} 
and its eigenvalues are called {\it multipliers} of system (\ref{system}). Thus, the monodromy matrix is the value at 
$t=1$   of the fundamental matrix $X(t,0)$ 
defined by the initial condition $X(0,0)=I$, 
and the multipliers are the roots of the equation
\begin{eqnarray*}
\det \left[ X(1,0) - \mu I \right] =0\,. 
\end{eqnarray*} 
Due to Theorem~2.3.1~\cite{Eastham} there exist the   open instability intervals. The   
 Assumption ISIN states that  there exists the nonempty   open instability interval $\Lambda \subset (0,\infty)$ for equation (\ref{Eastham1.2.1}).  
\smallskip

One can find in \cite{Eastham,Magnus_Winkler} the detailed description of functions $\alpha =\alpha (t) $ and $q=q(t) $ satisfying this condition.
For instance, in Theorem~4.4.1~\cite{Eastham} one can find asymptotic formula, which allows to estimate the length of the instability intervals 
of the equation obtained from (\ref{Eastham1.2.1}) by Liouville transformation. 
Then, according to the next lemma 
 one can find in 
the instability interval $\Lambda $  a 
number\,  $\lambda $\,  such  that a non-diagonal element of the monodromy matrix 
does not vanish. Moreover, this property is stable under small
perturbations of $\lambda $.

\begin{lemma} \mbox{\rm  (\cite{yagdjian_birk})}
Let  \,  $b(t)$\, be defined on \, ${\mathbb R}$\, non-constant,  positive, 
smooth function, which is  $1$-periodic.  Then there exists an open 
subset $\Lambda  ^0 \subset \Lambda  $
such that  $b_{21} \neq 0$ for all $\lambda  \in \Lambda ^0$.
\end{lemma}  
\medskip

Next we use the
 periodicity of \,$b = b(t)$\, and the eigenvalues \, $\mu_0>1$,\,  
$\mu_0^{-1}<1$
of the matrix 
$\,X_\lambda (1,0) \,$   
to construct solutions of (\ref{Eastham1.2.1}) with prescribed values on a discrete
set of time. The eigenvalues of matrix $X_\lambda (1,0)$ are $\mu_0$ and 
$\mu_0^{-1}$ with $ b_{11} + b_{22} = \mu_0 + \mu_0^{-1}$. 
Hence  $(b_{11} - \mu_0) + (b_{22} -\mu_0) = - \mu_0 + \mu_0^{-1}$
implies $|b_{11} - \mu_0| + |b_{22} - \mu_0| \geq 
|(b_{11} - \mu_0) + (b_{22} -\mu_0)| = |\mu_0 - \mu_0^{-1}| > 0$.
This leads to
\[
\max \{ |b_{11} - \mu_0| , |b_{22} - \mu_0|\} 
\geq \frac{1}{2} |\mu_0 - \mu_0^{-1}| > 0\,.
\]
Without loss of generality we can suppose
\[
|b_{11} - \mu_0| 
\geq \frac{1}{2} |\mu_0 - \mu_0^{-1}| > 0\,, \quad |b_{22} - \mu_0^{-1}| 
\geq \frac{1}{2} |\mu_0 - \mu_0^{-1}| > 0\,,
\]
because of $b_{11} - \mu_0 = -(b_{22} - \mu_0^{-1})$. Further, 
\[
1- \frac{b_{21}}{\mu_0^{-1}-b_{22}}\frac{b_{12}}{\mu_0-b_{11}}
= (\mu_0 - \mu_0^{-1})\frac{1}{b_{22} - \mu_0^{-1}}\not= 0 \,.
\]

\begin{lemma}\mbox{\rm (\cite{yagdjian_birk})}
\label{L3} 
Let $W = W(t)$, $V = V(t)$ be two  solutions of the equation   
\[
 w_{tt} + \left( \lambda    \alpha (t)-q(t) \right) w  = 0 
\]
 with the parameter $\lambda $ 
such that $b_{21} \not= 0$ and $b_{22}\not=\mu_0^{-1} $. Suppose then that $W = W(t)$   takes  the initial data
\[
 W(0) = 0\,,\quad
 W_t(0) =  1  \,,
\]
and that $V = V(t)$ takes the initial data 
\[
V(0) = 1\,,\quad
 V_t(0) =  0  \,.
\]
Then for every positive integer number \,$M \in {\mathbb N}$\,  one has
\begin{eqnarray*}
 W( M) 
& = &
 \frac{b_{21}}
{\mu_0-\mu_0^{-1}}( \mu_0^M - \mu_0^{-M} ) \,,\\ 
 V( M) 
& = &
 - \mu_0^M\frac{(b_{22}-\mu_0^{-1})}{(\mu_0-\mu_0^{-1})}
 +  \mu_0^{-M} \frac{b_{21}b_{12}}{(\mu_0-b_{11})(\mu_0-\mu_0^{-1})} .
\end{eqnarray*}
\end{lemma}

 For more applications of the Floquet-Lyapunov theory to hyperbolic equations with oscillating coefficients see \cite{Nishitani-Yagdjian,Reissig-Yagdjian, Ueda} and the bibliography therein. On the other hand, to study   the hyperbolic equations with oscillating coefficients one can   appeal to the so-called method of zone (see, e.g., \cite{Herrmann,Hirosawa,Wirth,yagdjian_book} and the  bibliography therein). 

\medskip

\section{Proof of Theorem~\ref{MTH}. Construction of blow--up solution to the scalar PDE}
\label{S4}

If condition (\ref{1.3-Nakanishi-Ohta_K}) of Theorem~\ref{MTH} does not hold, then (\ref{Not_N-O}) is true, that is,  $a >- \infty$ or $b < \infty$.
\smallskip

\noindent
If $u(t,x)$ is a solution
of (\ref{20}) and  takes initial values (\ref{Cdata})
then the function  (\ref{1.4-Nakanishi-Ohta-K})
solves the linear equation (\ref{lintrans}) 
and takes initial values
\begin{equation}
\label{Cdv-K}
v(0,x) =  \int_0^{u_0(x)}  \exp \left( \int_0^s f(r) dr \right) ds \,, \quad
v_{t}(0,x) =  u_1(x) \exp \left( \int_0^{u_0(x)}  f(r) dr \right) \,.
\end{equation}
Now let us choose initial data with the positive numbers $S>2n$ and $M$ which will be chosen later 
\begin{eqnarray*}
u_0(x)
& = &
\frac{1}{M^S}\chi  \left( \frac{x}{M^2}  \right) \in C_0^\infty ({\mathbb R}^n), \\
u_1(x)
& = &
\frac{A}{M^S}\chi \left( \frac{x}{M^2} \right)
\exp \left( -\int_0^{u_0(x)}  f(r) dr \right) \cos (x\cdot y) \in C_0^\infty({\mathbb R}^n)  ,
\end{eqnarray*}
where
$y \in {\mathbb R}^n$, $|y|^2 = \la$,  $\la$ is from the instability interval stated by ISIN,  while $\chi \in C_0^\infty ({\mb R}^n) $
is a non-negative cut-off function, $\chi (x)=1$ when $|x| \leq 1$. The
number
$A =\pm 1$, which is independent of the large
parameter $M \in {\mb N}$,
will be chosen later. Let $u=u(t,x)$ be a classical
solution of (\ref{20}) which takes these initial data. Then the function   $v(t,x) = G \big( u(t,x)\big)$
 solves  equation (\ref{lintrans}) and at $t=0$ takes values
\begin{eqnarray*}
v(0,x) & =  &
\int_0^{\frac{1}{M^S}\chi  \left( \frac{x}{M^2}  \right) }
\exp \left( \int_0^s f(r) dr \right) ds \in C_0^\infty ({\mathbb R}^n) \,,  \\
v_{t}(0,x) & =  &
  \frac{A}{M^S}\chi \left( \frac{x}{M^2} \right)
\cos \left( x\cdot y  \right) \in C_0^\infty ({\mathbb R}^n) \,.
\end{eqnarray*}
Let $W=W(t) $ be a solution given by Lemma~\ref{L3}. Consider the function
\[
V(t,x)\!  = \!\!\int_0^{\frac{1}{M^S}}\!\!\!
\exp   \left( \int_0^s \!\!\!f(r) dr \right) ds + W(t)\frac{b^{n/2}(t)}{b^{n/2}(0)}\frac{A}{M^S} \cos (x\cdot y)
\in C^\infty([0,\infty ] \times {\mathbb R}^n) \,.
\]
Function $V(t,x)$ solves equation (\ref{lintrans})
 while
\[
V(0,x)= \int_0^{\frac{1}{M^S}}
\exp \left( \int_0^s f(r) dr \right) ds , \quad
V_t(0,x) = \frac{A}{M^S} \cos (x\cdot y)  \,\, \mbox{\rm for all}
\,\, x \in {\mathbb R}^n\,.
\]
On the other hand for the function $v(t,x)$ we have
\[
v(0,x)  =
 \int_0^{\frac{1}{M^S}}
\exp \left( \int_0^s f(r) dr \right) ds , \quad
v_{t}(0,x)  = \frac{A}{M^S}\cos (x\cdot y) \,\, \mbox{\rm when} \,\,
 |x| \leq M^2.
\]
The finite propagation  speed in the Cauchy problem
(\ref{lintrans}),\,(\ref{Cdv-K}) implies
\[
V(t,x) = v(t,x) \quad {\rm in}  \quad
 \Pi_M:=[0,M]\times \{x \in {\mathbb R}^n ; |x| \leq M^{3/2} \}
\]
for large integer $M$. Hence
\[
v(t,x) =  \int_0^{\frac{1}{M^S}}
\exp \left( \int_0^s f(r) dr \right) ds + W(t)\frac{b^{n/2}(t)}{b^{n/2}(0)}\frac{A}{M^S}\cos (x\cdot y)
 \quad {\rm in}  \quad \Pi_M \,.
\]
In particular,
\begin{eqnarray*}
v(M,0)  
& = &
 \int_0^{\frac{1}{M^S}}
\exp \left( \int_0^s f(r) dr \right) ds
 +
\frac{A}{M^S}  \frac{b_{21}}
{\mu_0-\mu_0^{-1}}( \mu_0^M - \mu_0^{-M} )\,.
\end{eqnarray*}
Assume now that $b < \infty $. Then  the global existence of $u$ means
\begin{equation}
\label{21}
v (t,x) = \int_0^{u (t,x)} \exp \left( \int_0^s f(r) dr \right) ds < b
\quad \mbox{for all} \quad t \ge 0, \,x \in {\mathbb R}^n \,.
\end{equation}
We choose\, $A =1$,\,and  $S$  such that for
$M$ large enough one has (\ref{de}) for $u_0$, $u_1$.  On the other hand, there is a number $t(M)\in [0,M] $ such that   
$v(t(M),0) > b$. The last contradicts (\ref{21}).  
The case of $a>-\infty$ can be discussed in similar way. 
The theorem is proved.  \hfill $\fbox{}$
\bigskip

\section{\bf Proof of Theorem~\ref{T4}} 
\label{S5}

Assume that the problem has a global solution $(u^1(x,t),\ldots,u^m(x,t)) \in C^\infty  $ for every initial data 
$(u^1_\ell(x ),\ldots,u^m_\ell(x )) \in C^\infty  ({\mathbb R}^n)\times \ldots \times C^\infty  ({\mathbb R}^n)$, $\ell=0,1 $. We are going to prove that  
  the Nakanishi-Ohta condition (\ref{1.3-Nakanishi-Ohta_K}) is  fulfilled.   Consider the system of the equations (\ref{4}), 
where $\Gamma_{j,k}^i (u) $   are $C^\infty$ functions satisfying condition  (\ref{DL}) and  
\[
u^i(0,x) =a^iu_0(x), \quad u_t^i(0,x) = a^iu_1 (x)\,, \quad \quad i=1,\ldots m, \quad x
\in {\mathbb  R}^n\,.
\]
Consider also the Cauchy problem (\ref{20}) for the scalar   equation with 
the initial conditions  
\[
u (0,x) =u_0(x), \quad u_t (0,x) = u_1 (x)\,,  \quad x
\in {\mathbb  R}^n\,.
\]
Then 
the uniqueness and existence theorem  and condition (\ref{DL}) imply 
\[
u^1(t,x) =a^1u(t,x),\, u^2(t,x) =a^2u(t,x),\,\ldots ,\,u^m(t,x)=a^mu(t,x)    
\]
for all $x \in   {\mathbb  R}^n$, \,  $t \geq 0$. Thus we have obtained the existence of the global solution for the Cauchy problem   the nonlinear hyperbolic scalar equation (\ref{20}). 

Now we turn to the   scalar equation of (\ref{20}), 
where $f (u) $, $b(t)$  are $C^\infty$ functions and  $f (u) $ is from condition  (\ref{DL}). 
The Hopf-Cole-Nakanishi-Ohta transformation  converted equation  of (\ref{20})  into the linear wave equation for  $v$ defined by (\ref{lintrans}).  
Since $G \in C^2({\mathbb R})$ and $G'>0$, there exists the inverse $H$ (\ref{H}) of $G$, 
where $a$ and $b$ defined by (\ref{ab}). 
We choose initial data 
\[
 u _0(x )  =0\,, \qquad  u_1 (x ) = 1 \,,
\]
then 
\[
v(0)=0,\quad v_t(0)=1
\]
and 
\[
    v_{tt} - n\frac{\dot{b}(t) }{b(t)}v_t= 0\,.
\]
The explicit formula  for the solution $ v$ implies 
\[
 \int_0^{u (t )} \exp \left( \int_0^s f(r) dr \right) ds=v( t)=b^{-n}(0)\int_0^t b^{n}(\tau )\,d\tau \rightarrow \pm \infty \quad \mbox{\rm as} \quad t \rightarrow \pm \infty \,.
\]
Hence the condition (\ref{1.3-Nakanishi-Ohta_K}) is fulfilled.   The theorem is proved.  \hfill $\fbox{}$

\section{Proof of Lemma~\ref{L1}}
\label{S6}

In some chart the geodesic satisfy the system of equations
\[
\frac{d^2 u^i}{d s^2}  (s) + \sum_{j,k=1}^m\Gamma ^i _{j k }(u^1(s),\ldots,u^m (s) ) \frac{d  u^j  }{d s  }  (s) \frac{d  u^k }{d s  } (s)=0 \quad \mbox{\rm for all}\quad  i=1,\ldots, m\,.
\]
For the smooth geodesic lying in   the segment $I $ of  the straight line  \\
${\mathbb L}= \{(a_1t,\ldots,a_mt)\,|\, t \in {\mathbb R} \}$ of the Riemannian manifold $M$   we have\\ $u^1(s)=a_1 u(s)$, $\ldots\,, u^m(s)=a_m u(s)  $ for all $s  \in [c,d]$ and
\[
\left(  \frac{d  u   }{d s  }  (s)\right) ^2\sum_{j,k=1}^m\Gamma ^i _{j k }(a_1 u (s),\ldots,a_mu (s) ) a_ja_k=- a_i\frac{d^2 u }{d s^2}  (s)  \quad \mbox{\rm for all}\quad  s\in [c,d], 
\]
$i=1,\ldots, m$.\,  The  constant speed property of geodesics imply
\[
\left( \frac{d  u }{d s }  (s) \right)^2  \sum_{j,k=1}^m h_{kj}(a_1 u(s)  ,\ldots,a_m u(s)   )a_ja_k=constant\,.
\]
Consequently, the   function $  {d  u   (s) }/{d s  }  $ has no zeros and we can set 
\[
\widetilde{f} (s)=- \frac{d^2 u }{d s^2}  (s)   \left(  \frac{d  u   }{d s  }  (s)\right) ^{-2}
\quad {\rm and}\quad  f (u(s))=\widetilde{f} (s)\,,
\]  
since the function $u=u(s) $ has an inverse. On the other hand such geodesic covers the segment  $I \subseteq \mathbb L $ with  the parameter $t=u(s)$.
It follows (\ref{DL}). 
 
Conversely, suppose that  (\ref{DL}) holds. We can assume that $I= \{(a_1t,\ldots,a_mt)\,|\, t \in [1,b] \}$. Then for the point $(a_1,\ldots,a_m ) \in I$ we can solve the Cauchy problem for the scalar equation
\begin{equation}
\label{11}
\frac{d^2 u }{d s^2}  (s) + f(u(s))\left(  \frac{d  u   }{d s  }  (s)\right) ^2=0
\end{equation}
with the initial condition
\[
u   (0) =1\,,  \quad \frac{d  u }{d s }  (0)=\widehat{\xi },  \quad \mbox{\rm where}\quad
\widehat{\xi }^2=\left( \sum_{j,k=1}h_{kj}(a_1   ,\ldots,a_m   )a_ja_k \right)^{-1}\,.
\]
Further, since the point $ (a_1 u(s) ,\ldots,a_mu(s) )$ belongs to the segment $I$ for all sufficiently small $s$, the relation (\ref{11})   together with (\ref{DL}) imply 
\[
a_i\frac{d^2 u }{d s^2}  (s) + \left( \sum_{j,k=1}^m \Gamma_{j,k}^i  (a_1 u(s) ,\ldots,a_mu(s) ) a_ja_k\right)\left(  \frac{d  u   }{d s  }  (s)\right) ^2=0\,.
\]
Thus, $ (  u_1(s) ,\ldots, u_m(s) )= (a_1 u(s) ,\ldots,a_mu(s) )$  is a geodesic. The existence and uniqueness theorem for the system of ordinary differential equations  guarantees that   two geodesics with
a common point and equal tangent at that point must coincide. Hence, the geodesic covers the segment  $I \subseteq \mathbb L $. The lemma is proved. 
\hfill $\square$

\begin{remark}
The Poincar{\'e} half-plane model  (see,e.g.,\cite{Nishitani-Yagdjian}) possesses vertical half-lines which are distinguished geodesics. 
Another interesting example of a Lorentzian manifold that possesses  half-lines, which are distinguished geodesics is the Schwarzschild spacetime in the Eddington-Finkelstein
coordinates (see, e.g., \cite[Sec. 8.3]{Ohanian-R}).
\end{remark}

\section*{Acknowledgement}
The research of N.M.L.-R. was partially funded by the Japan Society for the Promotion of Science during June 19 --August 20 of 2018. Special thanks to Dr. Fumihiko Hirosawa for his help and for graciously hosting N.M.L.-R. at Yamaguchi University in summer of 2018.

\end{document}